\documentclass[]{interact}

\usepackage[english]{babel}
\usepackage[latin1]{inputenc}
\usepackage{epic}
\usepackage{amssymb}
\usepackage{amsmath}
\usepackage{calc}
\usepackage{graphicx}
\usepackage{amsthm}
\usepackage{multirow}
\usepackage{tikz}
\usepackage[ruled,vlined]{algorithm2e}
\usepackage{anysize}
\usepackage{epstopdf}
\usepackage[caption=false]{subfig}

\usepackage[longnamesfirst,sort]{natbib}
\bibpunct[, ]{(}{)}{;}{a}{,}{,}

\newcommand{\rv}[1]{#1}

\marginsize{2.5cm}{2.5cm}{2.5cm}{2.5cm} 

\begin{document}
	\title{Assessing the effectiveness of park-and-ride facilities on multimodal networks in smart cities}

	\author{
		\name{Juan A. Mesa \textsuperscript{a}\thanks{Juan A. Mesa. E-mail: jmesa@us.es}, Francisco A. Ortega \textsuperscript{b}\thanks{Francisco A. Ortega. E-mail: riejos@us.es}, Miguel A. Pozo \textsuperscript{a}\thanks{Miguel A. Pozo. E-mail: miguelpozo@us.es} and Ram\'on Piedra-de-la-Cuadra\textsuperscript{a}\thanks{Ram\'on Piedra-de-la-Cuadra. E-mail: rpiedra@us.es} }
		\affil{\textsuperscript{a}Departamento de Matem\'atica Aplicada II. Universidad de Sevilla, Spain; \textsuperscript{b}Departamento de Matem\'atica Aplicada I. Universidad de Sevilla, Spain}
	}
	
	\maketitle

	\begin{abstract}
		This paper presents an optimization procedure to choose a parking facility according to different criteria: total travel time including transfers, parking fee and a factor depending on the risk of not having an available spot in the park facility at the arrival time. An integer programming formulation is proposed to determine an optimal strategy of minimum cost considering the available information, different scenarios and each user profile. To evaluate the performance, a computational experience has been carried out on Seville (Spain), where an historical city center restricts the circulation of private vehicles and obliges the use of parking facilities.

	\end{abstract}


	\section{Introduction} 
	
	One of the key transport issues facing many countries in the world
	is the increase in congestion in urban areas and their accesses,
	due to the heavy dependence on the use of private cars. Since
	park-and-ride facilities allow citizens to use their cars for a
	part of the trip, while completing the rest by means of public
	transport, such facilities have been identified as a tool that can
	help to reduce traffic congestion \cite{Mes01}. During recent
	decades increased attention has been paid to park-and-ride
	facilities. Reports commissioned by urban or transportation
	agencies of cities, metropolitan areas or governments, and papers
	published in specialized journals, have informed about inquires and
	research related to several aspects of park-and-ride
	facilities. In the report \cite{Gle17} typical objectives of
	park-and-ride are listed.
	
	There are several variables related to the setting of
	park-and-ride facilities, being location, size and price, among which are the most
	important. In order to fix these variables the behavior decision
	of  potential users has been extensively studied. In fact,
	numerous studies meditate \cite{Bon04} on the complexity of users' preferences for choosing park-and-ride. Most parking choice models
	have applied logit models with data from
	stated choice surveys, conducted on resident and non-resident
	drivers. For example \cite{Hes04} have investigated parking choice
	based on a stated preference data set in the center of several cities in
	the United Kingdom using mixed logit models. \cite{Hab12} have
	investigated the relationship between parking choice and an activity
	scheduling process based on data coming from a survey in Montreal.
	Based on Decision Field Theory \cite{Qin13} analyze park-and-ride
	decision behavior in order to assist policy makers for planning
	park-and-ride facilities. An investigation on drivers' behavior has been
	carried out in a town in the north of Spain in \cite{Ibe14}. In
	that research three scenarios were considered: free parking
	on street, paid street parking and paid parking in an
	underground parking lot, and moreover three variables (access time to parking,
	parking fee and access time to destination) which were evaluated by
	means of a mixed logit model, where the values of time have been
	inferred. On the other hand, factors affecting mode change behavior of commuters were
	explored by using a multinomial logit in \cite{Isl15}.
	
	The evidence has shown that the attractiveness of a park facility
	is related to different cost attributes which can be combined with
	relative importance weights to build an overall measure of
	disutility. In \cite{Tho98} three types of costs were
	distinguished: 1) Access cost including in-vehicle travel time
	from the vehicle's current location to the car park as well as the
	time spent searching for an empty parking space, 2) waiting time cost
	that occurs when drivers have to make a line at a car park before being
	able to enter, 3) usage cost associated with a car park. When assessing
	these costs, travel time cost on foot or public transit to the
	destination has to be included. Moreover, the attractiveness of a
	parking lot depends on its condition, such as the size and the
	chance of finding a vacant parking space \cite{Hen01}. The relationship between
	prices and capacity of the parking lot has been studied in
	\cite{Gar02}, in which a bilevel programming model with the
	decisions on size and fees at the upper level and the reaction of
	the users at the lower level, has been proposed. The problem is solved
	by simulated annealing. Specifically, the number of available places 
	at each park facility for a given time could be estimated from the
	current occupancy, by considering the observed evolution of this
	parameter in similar settings of the reference day, taking the
	proximity to peak hours of traffic into account \cite{Hos14}.
	
	Searching-for-parking traffic comprises a significant amount of
	the total traffic volume, with average values reported equivalent
	to 30-50 $\%$  of the total traffic volume. Parking choice can be
	considered as a search process where drivers make a number of
	linked decisions based on updated knowledge gained by experience
	\cite{Tho98}.  According to studies
	carried out in several big cities, the daily parking-search time
	is evaluated at eight minutes on average \cite{Sho06}. This fact
	was validated from surveys conducted in Europe, which show that
	25-40 $\%$ of the travel time to central urban areas is attributed
	to the search for parking \cite{Bon04}. Hence, by reducing the
	parking-search time could lead to significant improvements in
	terms of travel time, production, traffic flows, fuel consumption,
	pollution and noise emission. This is one
	of the reasons why the searching for parking process has been
	investigated in recent years. \cite{Ben08} have developed an
	agent-based model to explicitly simulate the spatial search
	process and navigation behavior of drivers. In \cite{Kap12} a
	hierarchical structure of preferences to determine the chosen
	parking location, as well as the preferred route towards this
	destination has been presented. The study of
	drivers' behavior regarding uncertainties about search times for
	finding vacant spot has been carried out in \cite{Cha15}, where
	stated preference experiments are applied to several discrete
	choice models. One of the findings of this research is that most
	of the drivers searching for a parking lot who make a trip for
	shopping purposes start the searching process when they approach
	or reach their destinations. This fact supports the need for
	reducing the cruising traffic by introducing reservation systems
	and guiding tools. Cruising for parking time has been recently
	researched \cite{Arn17}. An evaluation of the parking search
	process when drivers are considered as agents of bounded
	rationality, assuming that their decisions are governed by
	lexicographic heuristics, is conducted in \cite{Kar17}.

	Users of parking lots can be divided into two groups: those using
	a subscription or a reservation with one  parking lot and non-reservation clients.
	The first ones pay a previously fixed discounted fees  but the
	second ones  pay according to the time they use the parking lot.
	In \cite{Gua09} a linear integer programming model for maximizing the revenue of a
	parking lot that considers several kinds of subscribers and individual users for
	different scenarios is analyzed. In this paper the forecast for
	each period of time is adjusted using a statistical methodology. A
	related problem that also takes into account the satisfaction of
	the drive-in customers is researched in \cite{Akh14}. The
	approach consists in combining optimization with simulation
	techniques.
	
	During the past two decades, traffic authorities of some cities
	and several companies have implemented devices to inform and guide
	drivers towards parking facilities with variables message signals
	in real time. These signaling systems show if park facilities are
	available and/or the number of parking places available (see Figure 1).

	\begin{center}
		\begin{figure}[h!]
			\centering
			\includegraphics[width=0.4\textwidth]{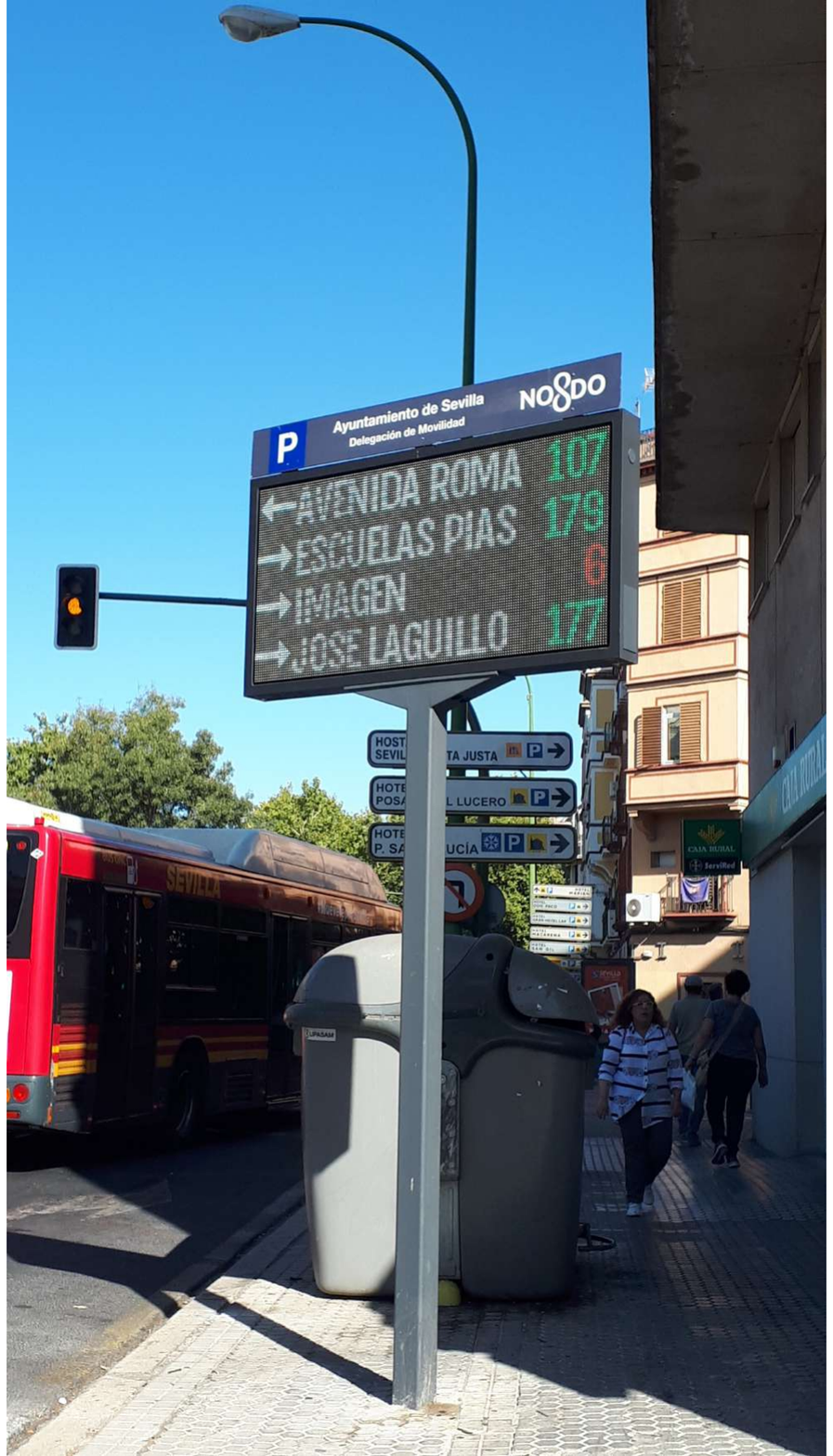}
			\caption{Signaling system of park facilities in the city of Seville}
			\label{fig:1}
		\end{figure}
	\end{center}

	Following
	these signals along the main roads, the drivers are guided along
	their routes towards the car parks. Moreover, Intelligent Transport
	Systems (ITS) have developed applications for helping drivers to be 
	informed and guided about the available parking facilities and,
	additionally, allowing them to make a reservation within the set of
	free spaces \cite{Tho97}. In \cite{Tho98} a simulation parking
	search model is introduced. Individual searching for a parking
	space accepts or rejects the vacant parking places based on a
	disutility function. In recent years, navigation systems of
	smart-phone applications have been trying to go beyond 
	suggesting the shortest or fastest paths in static scenarios, and have
	adopted algorithms to compute the shortest distance and fastest paths
	taking into account possible traffic jams. Note that if the
	traffic conditions were modified by a traffic accident or sudden
	congestion, the last user  request for searching for quickest  path
	should be re-calculated taking alternative routes into account.
	This aim is often accomplished by using historical data and
	pre-computed paths \cite{Far16}. Since the inspection time that the
	drivers spend while searching for an empty space increases when the
	number of available places is smaller, the expected waiting time
	depends on, among other factors, the occupancy level of the
	park facility. The decision of whether or not to accept the offered parking space will 
	affect the duration of a search. Therefore, this
	travel planning problem can be defined as a time-dependent
	shortest path problem through a fixed sequence of nodes
	\cite{Ber06}.

	
	In this paper, we consider a user that aims to travel from an origin to a destination splitting his trip into two parts. The first part will be carried out by private vehicle mode (car or motorcycle) ending in a park-and-ride facility to be determined. The second part of the trip, up to the final destination, will be carried out by means of public mode (including bus, tram, train or subway) and/or walking mode. Users differ from each other according to the fact of looking for a reservation in a parking facility and from the information they have. In this sense, users may be able to estimate the number of vacant places that parking facility will offer at a future time. Users may also check from a device the free places at a parking facility at a present time, which could be useful for estimating the availability of parking spaces at the arrival time at the parking facility. According to the type of user (and their available information), three criteria will have an influence on the decision towards choosing a park-and-ride facility: (1) the total travel time from origin to destination, (2) the parking fee and, (3) the attractiveness of the parking facility as a factor that will depend on the risk of not having an available spot at the parking facility at the arrival time.
	
	
	The remainder of this paper is structured as follows: Section 2 identifies some research lines devoted to network optimization related to this setting. Section 3, provides the description of the problem and a mathematical programming model that integrates the different aspects that are considered in the decision of the user. Computational results are presented in Section 4, followed by conclusions in Section 5.

	\section{Literature review}
	
	When planning a route in a road network, vehicle navigation systems usually provide the option of choosing among several different optimization criteria, or cost functions, such as a shortest route and a fastest route by giving preference to highways or avoiding toll highways.  In real-life scenarios the quality of a planned route is influenced by the presence of traffic jams along the suggested itinerary. Information on traffic jams, road works and incidences affecting road conditions can be received via smart-phone applications, so that the database, which describes the cost of crossing every arc in the road network, remains continuously updated in the GPS navigation system.
	Classical shortest path problems with fixed arc lengths have been intensively studied, generating several efficient algorithms (e.g, ~\cite{Dij59}, ~\cite{Dre69}, ~\cite{Ahu90}). In fastest path problems over a large dynamic network in real-time, the cost of an arc is the travel time of its arc and their objective is to find paths having minimum length with respect to time-dependent travel costs. 
	~\cite{Coo66} proposed a fastest path problem in networks where internode time requirements were included. 
	~\cite{Dre69} proposed a modification of the Dijkstra's static shortest path algorithm to calculate fastest paths between two nodes for a given departure time, by assuming that the travel times on the arcs are positive integers for every time period. 
	
	From a multiobjective perspective, ~\cite{Mar84} developed a multiple labeling version of  Dijkstra's label setting algorithm to generate all Pareto shortest paths from the source node to every other node in the network. ~\cite{Ham06} have generalized the classical shortest path problem by considering two objective functions in a setting of time-dependent data. In that paper a complete survey of the relevant literature is provided and an algorithm for the time-dependent bicriteria shortest path problem with non-negative data developed.
	
	The multimodal fastest path problem consists of finding a path from an origin to a destination while the total associated cost is minimized by means of the use of several transportation modes, such as personal car, taxi, subway, tram, bus, and walking.  ~\cite{Mod98} introduced an approach based on the classical shortest path problem for finding multi-objective shortest paths in urban multimodal transportation networks, taking the overall cost time and the users' inconveniences into consideration. In a multimodal network a node is a place where one has to select between continuing with the current mode or changing it. A change of mode or modal transfer is represented by an arc called transfer arc. An arc connecting two nodes by only one travel mode, is called travel arc. ~\cite{Loz01} have introduced label correcting techniques in an ad hoc multimodal shortest path algorithm to find the shortest viable path, where paths with illogical sequences in the use of transportation modes are eliminated. The same authors have extended their algorithm to calculate the viable hyperpath ~\cite{Loz02}.  The problem of finding the non-dominated viable shortest paths, considering the minimization of the travel time and of the number of modal transfers, is addressed in ~\cite{Art13}. 
	
	Finally, we must recall that the goal of this setting is to find a fastest path between two points in a road graph that additionally guarantees the visit of at least one intermediate node, which corresponds to a park facility. The Traveling Salesman Problem (TSP) consists of finding the minimum cost Hamiltonian tour on a given graph. This is one of the most studied combinatorial optimization problems. The Generalized Traveling Salesman Problem (GTSP) (see ~\cite{Lap83}) is a variation of the TSP in which the candidate node set to be visited is partitioned into clusters and the objective is to find a minimum cost Hamiltonian cycle passing through at least one node from each cluster. The GTSP and its variants arise in real-life applications such as computer operations, manufacturing logistics, distribution of goods by sea to the potential harbors and ecotourism routes (\cite{Lap96}, ~\cite{Bar16}). 
	Since the user must pass through an intermediate node (one parking garage between several possible options) before arriving at the destination point to leave his/her vehicle, we can assimilate this situation with the need of finding a minimum length route (not a tour) which includes exactly one visit to the set of parking lots. 
	
	In addition to these above mentioned research streams, several authors (see for instance ~\cite{Fli09}) have suggested improving the smart-phone applications for navigating in the city to the optimal park-and-ride facility by means of the incorporation of different aspects that have influence in the preferences of the drivers, such as parking costs, occupancy rates, walking time, possibility of shuttle buses, etc. To this end, this paper presents an optimization model that provides the best strategy for the user who requires guiding their navigation towards the most convenient park-and-ride facility in a metropolitan environment by considering the costs involved and the availability of empty places in the parking lot.


	\section{Problem description and formulation}
	
	\subsection{Problem description}
	
	We consider a user that aims to travel from an origin to a destination splitting his trip into two parts. 
	The first part will be carried out by private vehicle mode (car or motorcycle) ending in a park-and-ride facility to be determined. 
	The second part of the trip, up to the destination, will be carried out by means of public mode (including bus, tram, train or subway) and/or walking mode. 
	Once the park-and-ride facility has been chosen, the shortest time paths from origin to park-and-ride and from park-and-ride to destination, can be optimally determined.
	
	As previously mentioned, users differ from each other, firstly according to the fact of looking for (or not) a reservation at a parking facility. Secondly, the information that users could have may differ. In this sense, users may be able to estimate the number of free places that parking facility will offer at a future time according to the knowledge/experience that the user may have parking in a given area. Users may also check from a device the free places at a parking facility at a present time, that could be useful for estimating the availability of parking spaces at the arrival time at the parking facility. If such information could not be provided by the device, maybe it could show if the parking facility is full or not at a present moment. In any case, we assume that users always know the size of each parking facility.
	
	According to the type of user (and their available information), three criteria will have influence on the decision towards choosing a park-and-ride facility: (1) the total travel time from origin to destination, (2) the parking fee and (3) the ``attractiveness'' of the parking facility. We understand attractiveness as a factor (for users without a reservation in a parking facility) that will depend on the risk of not having an available spot at the parking facility at the arrival time. Other characteristics of the parking facility such as safety, parking space, cleaning, etc., could be taken into account. However, in order to keep the model as simple as possible, we have not included these features.
	
	
	
	Let us consider an initial directed graph $G= (N,A)$, where $N$ is the node set and $A$ is the arc set. We decompose $N=\{o\}\cup V\cup P\cup\{d\}$, where:
	\begin{itemize}
		\item $o$: node associated to the origin of the path, that is, the geographical position of the private vehicle mode (car or motorcycle) in the city.
		\item $d$:  node associated to the destination site in the city center, whose access is forbidden to cars.
		\item $P=P^+\cup P^-$ where  $P^+$ and $P^-$  are respectively  the entry nodes and the exit nodes of the car parks.
		\item $V$ is the set of intermediate nodes where the different transportation modes operate. 
	\end{itemize}
	
	In order to contruct the arc set $A$ we take into account the following issues:
	
	\begin{enumerate}
		\item There are no arcs entering at node $o$ and all arcs leaving this node correspond to the private vehicle mode.
		\item There are no arcs leaving node $d$ and all arcs entering in this node correspond to walking mode or other public mode.
		\item All arcs entering at a parking site $k\in P^+$ correspond to the private vehicle mode. From each node $k\in P^+$ only a single arc connects with the exit node $k'\in P^-$ of the parking site. Therefore, note that if $k\in P^+$ and $k'\in P^-$ are connected by an arc $(k,k')$ then both $k$ and $k'$ belong to the same park site. Additionally, all arcs leaving node $k'\in P^-$ correspond to walking mode or other public mode.
		\item All arcs entering and departing from a node $i\in V$ correspond to the same transportation mode. If there were a geographical node in which different transportation modes could be feasible for arriving at/departing from that node, it would be virtually replicated, once for each feasible transport mode. 
	\end{enumerate}
	
	Note that according to the graph construction, each feasible path connecting $o$ and $d$ must traverse a parking facility where the private vehicle mode changes to walking mode or other public mode. Additionally, several transportation modes could be included after leaving the parking facility just by means of including the corresponding transfer nodes and arcs between the different modes.
	
	
	In this paper, we assume an operational time horizon that is conveniently discretized in a set of time slots $T$. 
	Let $l_{ij} \in \mathbb{R}$ be the length associated with arc $(i, j)\in A$ and $v_{ij}^t$	the	average transit speed along the arc $(i,j)$ entering at node $i$ at time $t\in T$. For the case of $(k,k') \in A: k\in P^{^+},k'\in P^{^-}$, $l_{kk'}$ represents the transit time within parking facility $k$. 
	
	We define next some parameters related to parking facilities. Let $f_k^t$ be the fee at parking facility $k$ at time $t$ and $w_{kj}^t$ be the waiting time at parking facility exit node $k\in P^-$ during the transfer to node $j$ at time $t$. In case of requiring a later means of transport, we assume a fare $g_j$ has to be paid by the user.
	As previously mentioned, users differ from each other according to the information they have available (or not available). In this sense, let 
	$q_k^{1,t}$	be the estimation of free places that parking facility $k$ will offer at a future time $t$, this parameter is related to the knowledge/experience the user may have parking at $k$. 
	The user may also check a device that could provide the free places $q_k^{2}$ at parking facility $k$ at a present time  that could be useful for estimating the availability of parking plots at the arrival time at the parking facility. If this information was not provided by the device, maybe it could be shown if the parking facility is full or not. In this sense, let $q_k^{3}$	be the size of parking facility $k$ if it is free at a present time; 0 otherwise. In any case, we assume that users always know the size of each parking facility $k\in P^+$, namely $q_k^{4}$ as well as the size $Q$ of the biggest parking facility.
	Given $\gamma_0, \gamma_1, \gamma_2, \gamma_3, \gamma_4\in[0,1]$  weighting parameters holding $\sum_{i=1}^4\gamma_i=1$,  we can define the attractiveness of a parking facility as 
	
	$$\gamma_0\frac{\gamma_1 q_k^{1,t}+\gamma_2 q_k^2 +\gamma_3 q_k^3 +\gamma_4 q_k^4 }{Q}$$
	
	that its value is within $[0,1]$. \rv{Note that there are some relations among parameters $\gamma_i, i\in \{0,1,2,3,4\}$ according to the user profile. 
		First, users who look for a reservation (weighted as $\gamma_0=0$) are not affected by information of parking capacities ($\gamma_1=\gamma_2=\gamma_3=0$). 
		Second, available information of free places at parking facility $k$ at a present time (weighted with a $\gamma_2\neq 0$) includes the information of availability of the parking at a present time (that is, $\gamma_2\neq 0$ implies $\gamma_3=0$ and conversely, $\gamma_3\neq 0$, implies $\gamma_2=0$). 
		In addition, each user might impose a minimum level of attractiveness $\Lambda\in[0,1]$ accepted for the chosen parking. 
		The reader may observe that parameters $\gamma_0, \gamma_1, \gamma_2, \gamma_3, \gamma_4, \Lambda$ require to be properly weighted/calibrated according to the particular user's features and preferences. This might require an empiric procedure that is out of the scope of this research.}

	
	\subsection{Mathematical programming model}
	For each $t\in T$ and each arc $(i,j)$ let$x_{ij}^t\in\{0,1\}$ be the binary variable that takes the value 1, if the arc $(i,j)$ is traversed at time $t$, and 0 otherwise. We present next an integer programming formulation to determine a ``minimum generalized cost'' path from node $o$ to node $d$ in our setting. We understand that the path cost is determined by the minimization of three weighted criteria; namely, (1) the travel time, (2) the parking fee and (3) the non-attractiveness of the parking facility. For this aim, let $\alpha, \beta, \gamma_0$ be the weighting parameters for each of these criteria respectively.

	\hspace*{-2cm}\begin{eqnarray}
	\min && \alpha\sum_{t \in T}\left(\sum_{(i,j) \in A}  \frac{l_{ij}}{v_{ij}^t}\ x_{ij}^t +\sum_{(k,j) \in A: k\in P^{^-}}w_{kj}^t x_{kj}^t  \right)+\notag\\
	&&\hspace{0.5cm}+\beta\sum_{t \in T}\left(\sum_{(k,k') \in A: k\in P^{^+},k'\in P^{^-}}\hspace{-0.4cm}f_k^t x_{kk'}^t + \sum_{(k',j) \in A: k'\in P^{^-}}g_j x_{k'j}^t \right)+ \notag\\
	&&\hspace{0.5cm}+\gamma_0\left(\sum_{t \in T}\sum_{(i,k) \in A: k\in P^{^+}}\hspace{-0.6cm}x_{ik}^t
	\left(1-
	\frac{\gamma_1 q_k^{1,t}+\gamma_2 q_k^2 +\gamma_3 q_k^3 +\gamma_4 q_k^4 }{Q}
	\right)\right) \label{PAR1}\\                                                                
	\mbox{s.t.} && \sum_{t \in T}\sum_{(o,j) \in A} x_{oj}^t  =1,                  \label{PAR2}\\
	\mbox{  } && \sum_{t \in T}\sum_{(i,d) \in A } x_{id}^t  =1,   \label{PAR3}\\
	&& \sum_{(j,i) \in A} x_{ji}^t  -   \sum_{(i,j) \in A} x_{ij}^t  =0, \quad   i \in V \cup P, t \in T \label{PAR4}\\
	&& x_{ij}^t \in \{0,1\}, \quad     (i,j) \in A, t \in T. \label{PAR5}
	\end{eqnarray}
	
	The objective function \eqref{PAR1} minimizes a weighted sum of three criteria/components, namely travel time, parking cost and parking non-attractiveness.
	Constraint \eqref{PAR2} guarantees that the path will begin at the origin (by using the car as a means of transportation). 
	Constraint \eqref{PAR3} guarantees that the path will end at the destination (by using a means of an alternative transportation to the car).
	Constraints \eqref{PAR4} are flow conservation constraints for node subsets $V$ and $P$.
	
	We recall that according to the graph construction there is no need to impose on the model that no arcs enter in node $o$, no arcs leave from node $d$ and one parking facility has to be visited.
	
	Next, we propose two additional constraints to bound the travel time and the parking facility attractiveness: 
	
	\hspace*{0cm}\begin{eqnarray}
	\mbox{  } && \sum_{t \in T}\sum_{(i,d) \in A}  (t+\frac{l_{id}}{v_{id}^t}) x_{id}^t\leq \max\{t\in T\}, \quad  \label{PAR6}\\
	\mbox{  } && \sum_{t \in T}\sum_{(i,k) \in A: k\in P^{^+}}x_{ik}^t
	\frac{\gamma_1 q_k^{1,t}+\gamma_2 q_k^2 +\gamma_3 q_k^3 +\gamma_4 q_k^4 }{Q} 
	\geq \Lambda, \quad  \label{PAR7}
	\end{eqnarray}

	Constraint \eqref{PAR6} implies the arrival time to destination is upper bounded whereas 
	constraint \eqref{PAR7} ensures a minimal attractiveness for the chosen parking facility.

	If the arc costs induce no negative cycles on $G$, the problem of determining the best strategy for the user, who requires an optimal path between nodes $o$ and $d$ in this multimodal network, can be efficiently solved via ad-hoc polynomial time algorithms, like Dijkstra's algorithm.

	\section{Computational experiment}
	
	Next, we report on the results of a computational experiment that we have carried out in order to empirically show the model sensitivity by means of a parametric analysis of solutions.

	Seville is the capital of the Andalusia region (Spain), provided with a large and well-preserved historical center which is approximately 2,200 years old. In fact, the Historic Center of Seville is one of the largest in Europe, along with those of Venice and Genoa. It has an approximately circular configuration, with an area of 3.94 $km^2$. During 2016, the city of Seville received the visit of 2.5 million tourists attracted by an old town which contains three world heritage sites and also many convents, churches, palaces, museums and gardens. 
	Motorized traffic in this sector of Seville is limited and/or forbidden in many streets that are predominantly narrow and one-way.
	For this reason, there is a network of parking lots around the city center that provides, via signaling panels, information to the drivers about their current occupancy levels.
	We assume that this basic information to the users is complemented by the input data described along Section 3.
	
	\begin{center}
		\begin{figure}[h!]
			\centering
			\includegraphics[width=1.0\textwidth]{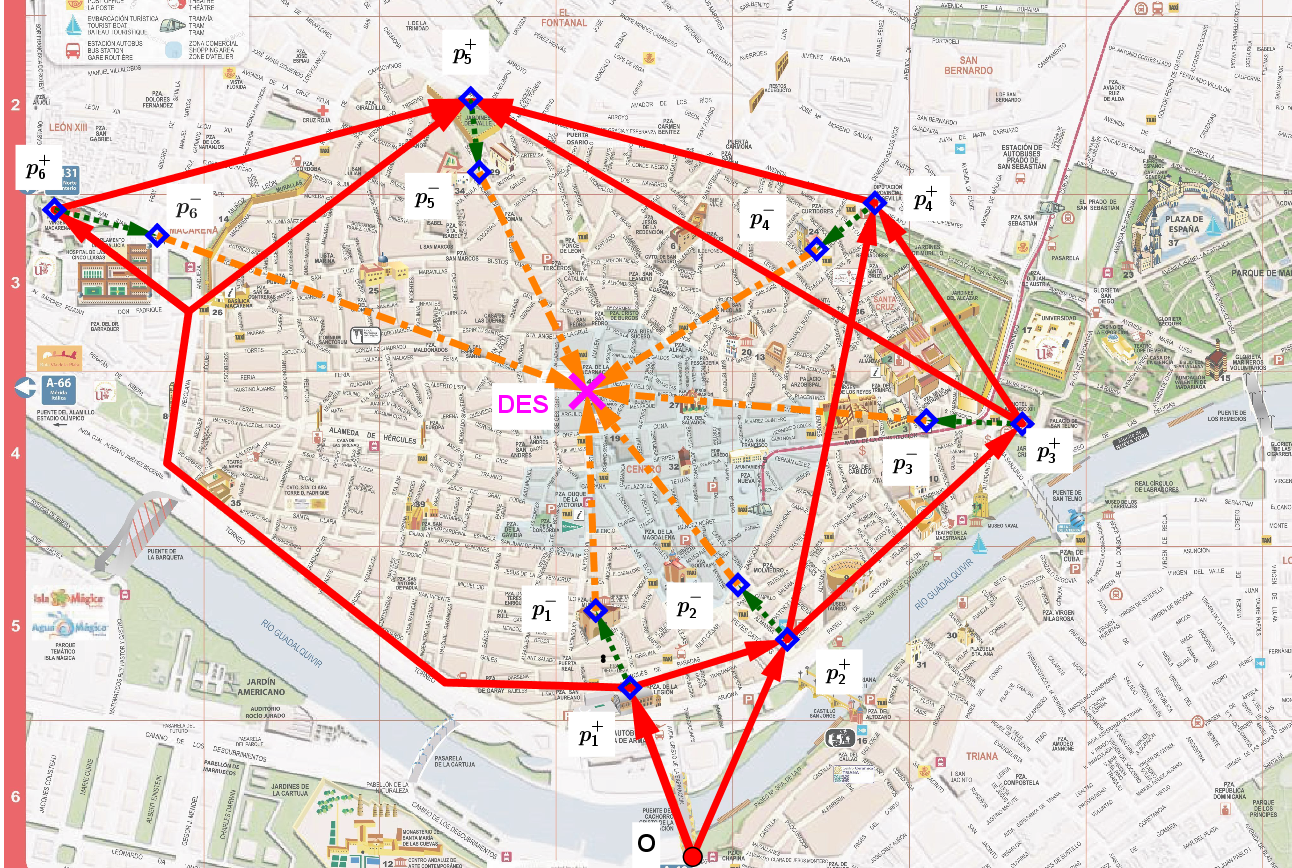}
			\caption{Graph associated to the selection of best strategy}
			\label{fig:2}
		\end{figure}
	\end{center}

	In Figure 2 we have supposed a driver with a starting point (labeled with $O$), and a destination (labeled with DES) located at the city center that is not accessible by private car. 
	Six existing car parks have been selected around the city center with entrances labeled by $p_k^{+}$ and exits $p_k^{-}$ for $k\in\{1,...,6\}$. Note that for the sake of visibility, exit parking nodes have been slightly separated from their real locations. 
	We deploy three different arc types according to the transportation mode in use: car (continuous red), walking (dotted green) and walking or bus (both in dashed orange).

	\begin{table}[h!]
		\tbl{Travel times involved in the fastest path from origin to destination.}
		{\begin{tabular}{|c|c|c|c|c|c|} \hline
				Park ($k$) & Time  to Park & Time in Park & Waiting time ($w_{kj}^t$) & Time to DES & Total time
				\\ \hline
				1       & 8  & 2 &  0 & 12  & 22
				\\ \hline
				2       & 7  & 1.5 & 0  & 13  & 20.5
				\\ \hline
				\textbf{3}  & 11 & 1 & 2  & 6  & \textbf{20}
				\\ \hline
				4       & 14 & 2 & 5  & 13  & 34
				\\ \hline
				5       & 15 & 1 & 2 & 9  & 27
				\\ \hline
				6       & 10 & 1 &  3.5 & 15  & 29.5
				\\ \hline
		\end{tabular}}
		\label{scenario_1}
	\end{table}
	
	Table \ref{scenario_1} shows traveling times in minutes from the origin point (La Pa\~noleta, 41910-Seville, Spain) to the destination point (Plaza del Salvador, 41004-Seville, Spain) passing through each one of the 6 available parking lots. 
	Column ``Time to Park'' shows the driving time from the origin to each parking entrance through a time dependent shortest-path. 
	Column ``Time in Park'' shows the time spent at the parking facility. 
	Column ``Waiting Time'' shows the waiting time at the exit of the parking facility in case a bus transfer is required. 
	Finally, column ``Time to DES'' shows the total travel time from the exit parking node to DES through a time dependent shortest-path that might include bus and/or walking mode.  
	Note that the best parking facility (marked in bold) in terms of total time is parking lot 3.
	In addition, Table \ref{scenario_1} compares travel times obtained when using bus mode vs walking mode after parking the car. Depending on the frequency of bus services, parking lots 1 or 2 might also be a good option.

	\begin{table}[h!]
		\tbl{Park fees and fare rates involved in the trip.}
		{\begin{tabular}{|c|c|c|c|} \hline
				Park ($k$)  & Parking fee ($f_k^t$)& Fare rate ($g_j$)& Total price 
				\\ \hline
				1       & 1.17  & 0 & 1.17 
				\\ \hline
				2       & 1.01  & 0 &  1.01  
				\\ \hline
				3       & 0.3 & 0.74 & 1.04  
				\\ \hline
				4       & 0.5 & 0.69 & 1.19  
				\\ \hline
				5       & 0.35 & 0.74 & 1.09 
				\\ \hline
				\textbf{6}       & 0.16 & 0.69 &  \textbf{0.85} 
				\\ \hline
		\end{tabular}}
		
		\label{t:price}
	\end{table}
	
	Table \ref{t:price} includes parking fees and fare rates if a later means of transportation is required. 
	Note that the best option (cheapest) in this case (parking 6) changes with respect to  best option (parking lot 3) in the previous case (fastest).
	
	\begin{table}[h!]
		\tbl{Parking parameters involved in the non-attractiveness criterion.}
		{
			\begin{tabular}{|c|c|c|c|c|c|       } \hline
				Park ($k$)  & Estimation of free places ($q^{1,t}_k$) & Free places ($q^2_k$) & Full or not ($q^3_k$) & Parking size ($q^4_k$) & $Q$
				\\ \hline
				1       & 15  & 0 &  0 & 25  & 100
				\\ \hline
				2       & 0  & 8 & 30  & 30  & 100
				\\ \hline
				3       & 20 & 1 & 50  & 50  & 100
				\\ \hline
				4       & 6 & 5 & 25  & 25  & 100
				\\ \hline
				5       & 20 & 50 & 70 & 70  & 100
				\\ \hline
				6       & 66 & 35 & 100 & 100  & 100
				\\ \hline
			\end{tabular}
		}
		
		\label{scenario_3}
	\end{table}

	Table \ref{scenario_3} shows the different parking parameters involved in the non-attractiveness criterion of the objective function. In order to obtain a non-attractiveness value, these parameters require to be combined with weights $\gamma_i$, $i\in \{0,1,2,3,4\}$ according to the individual specific preferences and characteristics of the user in consideration.

	\begin{table}[h!]
		\tbl{Weighting parameters and generalized cost for a user looking for a reservation.}
		{%
			\begin{tabular}{|c|c|c|c|c|c|c|c|c|c|c|c|      } \hline
				Park ($k$)  & $\alpha$ & Total time & $\beta$ & Total price &$ \gamma_0$ & $\gamma_1$ &$ \gamma_2$ & $\gamma_3$ & $ \gamma_4$  & Non-attractiveness & Total cost
				\\ \hline
				1    & 0.2  & 4.4 & 0.8 & 0.93  & 0 & 0 & 0 & 0 &1 & 0 & 5.33
				\\ \hline
				2    & 0.2  & 4.1 & 0.8  & 0.81   & 0 & 0 & 0 & 0 & 1 & 0 &4.91
				\\ \hline
				\textbf{3}    &0.2 & 4 & 0.8  & 0.83   & 0 & 0 & 0 & 0 & 1 & 0&\textbf{4.83}
				\\ \hline
				4    & 0.2 & 6.8 & 0.8  & 0.95   & 0 & 0 & 0 & 0 & 1 & 0&7.75
				\\ \hline
				5    & 0.2 & 5.4 & 0.8 &0.87   & 0 & 0 & 0 & 0 & 1 & 0&6.27
				\\ \hline
				6    &0.2 & 5.9 & 0.8 & 0.68   & 0 & 0 & 0 & 0 & 1 & 0&6.58
				\\ \hline
			\end{tabular}
		}

		\label{scenario_4}
	\end{table}

	\begin{table}[h!]
		\tbl{Weighting parameters and generalized cost for a user with information about free parking plots.}
		{%
			\begin{tabular}{|c|c|c|c|c|c|c|c|c|c|c|c|      } \hline
				Park ($k$)  & $\alpha$ & Total time & $\beta$ & Total price &$ \gamma_0$ & $\gamma_1$ &$ \gamma_2$ & $\gamma_3$ & $ \gamma_4$  & Non-attractiveness & Total cost
				\\ \hline
				1    & 0.04  & 0.88 & 0.17 &0.2  & 0.79 & 0 &0.67 & 0 & 0.33 & 0.72 & 1.80
				\\ \hline
				2    &  0.04  & 0.82 & 0.17  & 0.17   & 0.79 & 0 & 0.67 & 0 & 0.33 & 0.67 &1.66
				\\ \hline
				3    &  0.04 & 0.8 & 0.17  & 0.18   & 0.79 & 0 & 0.67 & 0 & 0.33 & 0.65&1.63
				\\ \hline
				4    &  0.04 & 1.36 & 0.17  &0.2   & 0.79 & 0 & 0.67 & 0 & 0.33 & 0.7&2.26
				\\ \hline
				\textbf{5}    &  0.04 & 1.08 & 0.17 &0.19   & 0.79 & 0 & 0.67 & 0 & 0.33 & 0.34&\textbf{1.61}
				\\ \hline
				6    &  0.04 &1.18 & 0.17 &0.14   &0.79 & 0 &0.67 & 0 & 0.33 &0.34&1.67
				\\ \hline
			\end{tabular}
		}

		\label{scenario_5}
	\end{table}
	
	Tables \ref{scenario_4} and \ref{scenario_5} provide the objective function values (generalized total cost) when the three criteria are combined according to the objective function \eqref{PAR1}. 
	In both cases weighting values $\alpha, \beta, \gamma_i, \, i\in \{0,1,2,3,4\}$ have been chosen for modeling to different user's profiles.
	In particular, Table \ref{scenario_4} shows weights of a user that chose a parking lot with a reservation procedure. For that reason, the objective value is not affected by information of parking capacities.
	Table \ref{scenario_5} shows weights of a user with knowledge of the number of free plots at the initial time as well as the size of the parking. Besides, this example assumes that this user is not able to estimate the number of free plots at his/her arriving time at the parking.

	If we take into account traffic congestion effects in the Seville metropolitan area at rush hour (see Figure 3) and its impact on the access time to parking lots 1 to 6, the arc cost must be recalculated in a part of the graph, following the approach of ~\cite{Pal03}. In the example being analyzed, the access times would increase by 3 minutes to reach the car parks 1 and 2, as well as 1 minute to access the car parks 5 and 6. These specific modifications in the travel times of certain arcs do not change the final decisions in the scenarios analyzed, although they could do so if the levels of congestion were more pronounced.
	
	\begin{center}
		\begin{figure}[h!]
			\centering
			\includegraphics[width=0.75\textwidth]{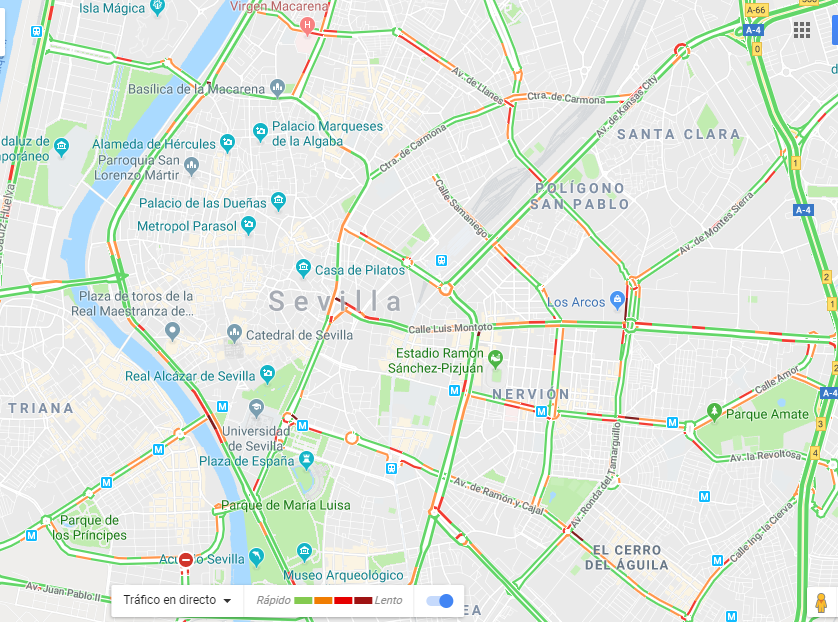}
			\caption{Map of traffic intensities in the city of Seville}
			\label{fig:3}
		\end{figure}
	\end{center}

	
	\section{Conclusions}
	Intelligent car parking systems, through smart-phone applications, are nowadays evolving to provide the best park facility for a given destination in terms of distance and occupancy level. In order to improve the navigation provided by these applications, a global assessment of different aspects that have influence in the drivers' preferences (like parking costs, walking time, vacant places forecast, transfer possibilities with other means of transportation, etc.) can be integrated when a user needs to necessarily visit a park-and-ride facility before reaching his/her destination located in a restricted traffic zone.
	
	Determining the best strategy for a user, who requires navigating towards the best choice of park-and-ride facility in an urban setting, has been dealt with in this paper. The proposed approach has consisted of formulating an integer programming model where a generalized cost is minimized, while a constraint set guarantees feasibility of solutions into the ad hoc built graph. The optimization model can be solved by using an adapted shortest path algorithm. 
	
	The main contribution of the paper is twofold. First, the objective in the model integrates several cost attributes: determining fastest routes through a multimodal network with time-dependent transit times, minimizing parking costs and an attractiveness criterion related to the risk of not having an available space at the chosen parking facility.
	Second, the coefficient values allow us to collect, in the same optimization model, the diversity of several user profiles, as well as the information available for the different scenarios.
	Results of a computational experiment located in the urban setting of Seville (Spain), carried out in order to empirically show the model sensitivity to the input parameters, have been reported.
	

	\section*{Acknowledgements}
	This research has been partially supported by the Spanish Ministry of Economy and Competitiveness through grants MTM2015-67706-P (MINECO/FEDER, UE) and MTM2016-74983-C2-1-R (MINECO/FEDER, UE). This support is gratefully acknowledged.

\end{document}